\documentclass[10pt]{amsart}
\usepackage{amssymb}
\usepackage{amsthm,amsmath}
\usepackage{cite}

\usepackage{ccfonts}

\usepackage[T1]{fontenc}

\title[A Generalization of Bohr-Mollerup's Theorem]{A Generalization of Bohr-Mollerup's Theorem For Higher Order Convex Functions: A~Tutorial}

\author{Jean-Luc Marichal}\thanks{Corresponding author: Jean-Luc Marichal, University of Luxembourg, Department of Mathematics, Maison du Nombre, 6, avenue de la Fonte, L-4364 Esch-sur-Alzette, Luxembourg. Email: jean-luc.marichal@uni.lu}
\address{University of Luxembourg, Department of Mathematics, Maison du Nombre, 6, avenue de la Fonte, L-4364 Esch-sur-Alzette, Luxembourg}
\email{jean-luc.marichal[at]uni.lu}

\author{Na\"im Zena\"idi}
\address{University of Li\`ege, Department of Mathematics, All\'ee de la D\'ecouverte, 12 - B37, B-4000 Li\`ege, Belgium}
\email{nzenaidi[at]uliege.be}

\date{\today}

\begin{document}

\theoremstyle{plain}
\newtheorem{theorem}{Theorem}[section]
\newtheorem{lemma}[theorem]{Lemma}
\newtheorem{proposition}[theorem]{Proposition}
\newtheorem{corollary}[theorem]{Corollary}
\newtheorem{fact}[theorem]{Fact}

\theoremstyle{definition}
\newtheorem{definition}[theorem]{Definition}

\newtheorem{ex}[theorem]{Example} 
    \newenvironment{example}    
    {\renewcommand{\qedsymbol}{$\lozenge$}%
    \pushQED{\qed}\begin{ex}}
    {\popQED\end{ex}}

\theoremstyle{remark}
\newtheorem{rem}[theorem]{Remark} 
    \newenvironment{remark}    
    {\renewcommand{\qedsymbol}{$\lozenge$}%
    \pushQED{\qed}\begin{rem}}
    {\popQED\end{rem}}

\theoremstyle{remark}
\newtheorem*{claim}{Claim}
\newtheorem*{conjecture}{Conjecture}

\newcommand{\R}{\mathbb{R}}
\newcommand{\N}{\mathbb{N}}
\newcommand{\Z}{\mathbb{Z}}
\renewcommand{\S}{\mathrm{S}}
\newcommand{\cC}{\mathcal{C}}
\newcommand{\cD}{\mathcal{D}}
\newcommand{\cK}{\mathcal{K}}
\newcommand{\cR}{\mathcal{R}}
\newcommand{\ran}{\mathrm{ran}}
\newcommand{\emptybox}{\mathstrut^{\mathstrut}_{\mathstrut}}
\newcommand{\li}{\,\mathrm{li}}

\def\tchoose#1#2{{\textstyle{{{#1}\choose{#2}}}}}
\def\parag#1{\medskip\noindent{\textbf{#1.}}}

\begin{abstract}
In its additive version, Bohr-Mollerup's remarkable theorem states that the unique (up to an additive constant) convex solution $f(x)$ to the equation $\Delta f(x)=\ln x$ on the open half-line $(0,\infty)$ is the log-gamma function $f(x)=\ln\Gamma(x)$, where $\Delta$ denotes the classical difference operator and $\Gamma(x)$ denotes the Euler gamma function. In a recently published open access book, the authors provided and illustrated a far-reaching generalization of Bohr-Mollerup's theorem by considering the functional equation $\Delta f(x)=g(x)$, where $g$ can be chosen from a wide and rich class of functions that have convexity or concavity properties of any order. They also showed that the solutions $f(x)$ arising from this generalization satisfy counterparts of many properties of the log-gamma function (or equivalently, the gamma function), including analogues of Bohr-Mollerup's theorem itself, Burnside's formula, Euler's infinite product, Euler's reflection formula, Gauss' limit, Gauss' multiplication formula, Gautschi's inequality, Legendre's duplication formula, Raabe's formula, Stirling's formula, Wallis's product formula, Weierstrass' infinite product, and Wendel's inequality for the gamma function. In this paper, we review the main results of this new and intriguing theory and provide an illustrative application.
\end{abstract}

\keywords{Difference equation, higher order convexity, Bohr-Mollerup's theorem, principal indefinite sum, Gauss' limit, Euler product form, Raabe's formula, Binet's function, Stirling's formula, Gauss' multiplication formula, Euler's constant, gamma and polygamma functions.}

\subjclass[2010]{26A51, 33B15, 33B20, 39A06, 39B22.}

\maketitle

\section{Introduction}

One of the best-known special functions of mathematical analysis is the Euler gamma function. Its restriction to the real open half-line $\R_+=(0,\infty)$ is usually defined as the following improper integral (see, e.g., Srinivasan~\cite{Sri07})
$$
\Gamma(x) ~=~ \int_0^{\infty}t^{x-1}{\,}e^{-t}{\,}dt,\qquad x>0.
$$

It is well known and easily seen that this function satisfies $\Gamma(1)=1$ and the identity (using integration by parts)
$$
\Gamma(x+1) ~=~ x\,\Gamma(x),\qquad x>0.
$$

In 1922, Bohr and Mollerup~\cite{BohMol22} established the following simple, but remarkable characterization of the gamma function.

\begin{theorem}[Bohr-Mollerup's theorem]\label{thm:B-Madd44760}
The gamma function is the unique logarithmically convex solution $f\colon\R_+\to\R_+$ satisfying $f(1)=1$ to the equation \begin{equation}\label{eq:FunEq639}
f(x+1) ~=~ x{\,}f(x),\qquad x>0.
\end{equation}
\end{theorem}

A decade later, Artin \cite{Art31} (see also Artin \cite{Art15}) investigated and simplified the proof of this characterization, which has become also known as the Bohr-Mollerup-Artin theorem. He also showed that many classical properties of the gamma function are actually very elementary consequences of this theorem and its proof. Among these properties, recall the \emph{Stirling formula}
\begin{equation}\label{eq:Stir46}
\lim_{x\to\infty}\frac{\Gamma(x)}{\sqrt{2\pi}{\,}e^{-x}{\,}x^{x-\frac{1}{2}}} ~=~ 1,
\end{equation}
the \emph{Gauss multiplication formula}
\begin{equation}\label{eq:GaussMult}
\prod_{j=0}^{m-1}\Gamma\left(\frac{x+j}{m}\right) ~=~ \frac{\Gamma(x)}{m^{x-\frac{1}{2}}}{\,}(2\pi)^{\frac{m-1}{2}},\qquad x>0,~m=1,2,\ldots,
\end{equation}
and the \emph{Gauss limit}
\begin{equation}\label{eq:GaussLim6689}
\Gamma(x) ~=~ \lim_{n\to\infty}\frac{n!{\,}n^x}{x(x+1)\cdots (x+n)}{\,},\qquad x>0.
\end{equation}

It is not difficult to see that Bohr-Mollerup's theorem above can be slightly generalized as follows.

\begin{theorem}[Bohr-Mollerup's theorem]
All logarithmically convex solutions $f\colon\R_+\to\R_+$ to equation \eqref{eq:FunEq639} are of the form $f(x)=c{\,}\Gamma(x)$, where $c>0$.
\end{theorem}

Indeed, if $f\colon\R_+\to\R_+$ is a logarithmically convex solution to equation \eqref{eq:FunEq639}, then clearly so is the function $f/f(1)$, which must be the gamma function by Bohr-Mollerup's Theorem~\ref{thm:B-Madd44760}.

The following theorem provides a reformulation of the latter result using the additive notation, where $\Delta$ stands for the classical difference operator.

\begin{theorem}[Additive version of Bohr-Mollerup's theorem]\label{thm:B-Madd4476}
All convex solutions $f\colon\R_+\to\R$ to the equation $\Delta f(x)=\ln x$ are of the form $f(x)=c+\ln\Gamma(x)$, where $c\in\R$.
\end{theorem}

It is natural to ask whether analogues of Theorem~\ref{thm:B-Madd4476} can be obtained by replacing the logarithm function in the difference equation $\Delta f(x)=\ln x$ with any other real function. In a recently published monograph \cite{MarZen22}, the authors showed that such analogues do exist for a very wide variety of functions. We now state this result in the following uniqueness theorem, which actually constitutes a major generalization of Bohr-Mollerup's theorem.

Recall first that a function $f\colon I\to\R$, where $I$ is any nontrivial real interval, is said to be \emph{$p$-convex} (resp.\ \emph{$p$-concave}) for some integer $p\geq 0$ if for any pairwise distinct points $x_0,x_1,\ldots,x_{p+1}$ in $I$ we have that
$$
f[x_0,x_1,\ldots,x_{p+1}]~\geq ~0\qquad (\text{resp.}~f[x_0,x_1,\ldots,x_{p+1}]~\leq ~0),
$$
where the symbol $f[x_0,x_1,\ldots,x_{p+1}]$ stands for the divided difference of $f$ at the points $x_0,x_1,\ldots,x_{p+1}$. It can be shown that, if $I$ is an open interval and $f$ is $p$ times differentiable, then it is $p$-convex (resp.\ $p$-concave) if and only if $f^{(p)}$ is increasing (resp.\ decreasing). For background and references, see, e.g., \cite[Section 2.2]{MarZen22}.

We say that a function $f\colon\R_+\to\R$ is \emph{eventually $p$-convex} (resp.\ \emph{eventually $p$-concave}) if it is $p$-convex (resp.\ $p$-concave) in some neighborhood of infinity.

For any integers $p\geq 0$ and $n\geq 1$, and any function $g\colon\R_+\to\R$, we define the function $f^p_n[g]\colon\R_+\to\R$ by the equation
$$
f^p_n[g](x) ~=\null \sum_{k=1}^{n-1}g(k)-\sum_{k=0}^{n-1}g(x+k)+\sum_{j=1}^p\tchoose{x}{j}\,\Delta^{j-1}g(n),\qquad x>0.
$$

\begin{theorem}[Uniqueness]\label{thm:int1}
Let $p\geq 0$ be an integer and let the function $g\colon\R_+\to\R$ have the property that the sequence $n\mapsto\Delta^pg(n)$ converges to zero. Suppose that $f\colon\R_+\to\R$ is an eventually $p$-convex or eventually $p$-concave function satisfying the difference equation $\Delta f=g$. Then $f$ is uniquely determined (up to an additive constant) by $g$ through the equation
$$
f(x) ~=~ f(1)+\lim_{n\to\infty}f^p_n[g](x),\qquad x>0,
$$
and the convergence is uniform on any bounded subset of $\R_+$.
\end{theorem}

Taking $p=1$ and $g(x)=\ln x$ in Theorem~\ref{thm:int1}, we immediately retrieve both Bohr-Mollerup's Theorem~\ref{thm:B-Madd4476} and Gauss' limit \eqref{eq:GaussLim6689}. We thus see that Theorem~\ref{thm:int1} provides a generalization of Bohr-Mollerup's theorem to a vast spectrum of functions. The following example, which will be our guiding example throughout this paper, provides another illustration of this generalization.

\begin{example}[The polygamma function $\psi_{-2}$]\label{ex:PolyG}
Consider the polygamma function $\psi_{-2}\colon\R_+\to\R$ defined by the following equation (see, e.g., Adamchik \cite{Ada98})
$$
\psi_{-2}(x) ~=~ \int_0^x\ln\Gamma(t){\,}dt,\qquad x>0.
$$
It is known (see, e.g., Adamchik \cite[p.~196]{Ada98} and Remmert \cite[p.~46]{Rem98}) that this function satisfies
$$
\psi_{-2}(1) ~=~ \frac{1}{2}\ln(2\pi).
$$
Moreover, it is $2$-convex since its second derivative is the digamma function $\psi$, which is increasing on $\R_+$ (see, e.g., Srivastava and Choi \cite{SriCho12}). Furthermore, for any $x>0$ we have
\begin{eqnarray*}
\Delta f(x) &=& \int_0^1\ln\Gamma(t){\,}dt+\int_1^{x+1}\ln\Gamma(t){\,}dt-\int_0^x\ln\Gamma(t){\,}dt\\
&=& \psi_{-2}(1)+\int_0^x(\ln\Gamma(t+1)-\ln\Gamma(t)){\,}dt\\
&=& \psi_{-2}(1)+\int_0^x\ln t{\,}dt.
\end{eqnarray*}
Thus, we have $\Delta f=g$ on $\R_+$, where the function $g\colon\R_+\to\R$ is defined by the equation
$$
g(x) ~=~ x\ln x -x+\frac{1}{2}\ln(2\pi),\qquad x>0,
$$
and has the property that the sequence $n\mapsto\Delta^2g(n)$ converges to zero. It follows from Theorem~\ref{thm:int1} that the function $\psi_{-2}$ is the unique (up to an additive constant) eventually $2$-convex solution to the equation $\Delta f=g$ on $\R_+$.
\end{example}


In \cite[Chapter 3]{MarZen22} we also provided the following result, which gives sufficient conditions on the function $g$ for the existence of an eventually $p$-convex or eventually $p$-concave solution to the difference equation $\Delta f=g$.

\begin{theorem}[Existence]\label{thm:int2}
Let $p\geq 0$ be an integer and suppose that the function $g\colon\R_+\to\R$ is eventually $p$-convex or eventually $p$-concave and has the asymptotic property that the sequence $n\mapsto\Delta^pg(n)$ converges to zero. Then there exists a unique (up to an additive constant) eventually $p$-convex or eventually $p$-concave solution $f\colon\R_+\to\R$ to the difference equation $\Delta f=g$. Moreover,
\begin{equation}\label{eq:int2}
f(x) ~=~ f(1)+\lim_{n\to\infty}f^p_n[g](x),\qquad x>0,
\end{equation}
and $f$ is $p$-convex (resp.\ $p$-concave) on any unbounded subinterval of\/ $\R_+$ on which $g$ is $p$-concave (resp.\ $p$-convex). Furthermore, the convergence in \eqref{eq:int2} is uniform on any bounded subset of $\R_+$.
\end{theorem}

We observe that Theorem~\ref{thm:int2} was first proved in the special case when $p=0$ by John~\cite{Joh39}. It was also established in the case when $p=1$ by Krull \cite{Kru48,Kru49} and then in its multiplicative version by Webster \cite{Web97a,Web97b}.

We also observe that identity \eqref{eq:int2} actually provides for the function $f$ an analogue of Gauss' limit \eqref{eq:GaussLim6689} for the gamma function. More generally, for every solution $f\colon\R_+\to\R$ arising from Theorem~\ref{thm:int2}, we presented in \cite{MarZen22} counterparts of various classical properties of the gamma function, including analogues of \emph{Burnside's formula}, \emph{Euler's infinite product}, \emph{Euler's reflection formula}, \emph{Gauss' limit}, \emph{Gauss' multiplication formula}, \emph{Gautschi's inequality}, \emph{Legendre's duplication formula}, \emph{Raabe's formula}, \emph{Stirling's formula}, \emph{Wallis's product formula}, \emph{Weierstrass' infinite product}, and \emph{Wendel's inequality} for the gamma function. We also introduced and discussed analogues of \emph{Binet's function}, \emph{Euler's constant}, \emph{Fontana-Mascheroni's series}, \emph{Stirling's constant}, \emph{Webster's inequality}, and \emph{Webster's functional equation}. We also provided and discussed some additional properties, including asymptotic equivalences, asymptotic expansion formulas, Taylor series expansion formulas, and Gregory formula-based series representations.

All these properties, combined with the uniqueness and existence theorems above, actually offer a unifying setting that enables us to systematically investigate a very wide variety of functions. This fact was largely discussed in \cite{MarZen22} and even extensively illustrated through various examples, ranging from the gamma function itself and its best-known variants to important special functions such as the Hurwitz zeta function and the generalized Stieltjes constants.

In the present paper, we provide a summary of the main results of this new and intriguing theory. We also illustrate these results by applying them to the polygamma function $\psi_{-2}$ (see Example~\ref{ex:PolyG}), which will be our guiding example throughout.

The outline of this paper is as follows. In Section 2, we introduce the concept of the \emph{principal indefinite sum} from the solutions arising from the existence Theorem~\ref{thm:int2}. In Sections 3 to 9, we introduce the analogues of Gauss' limit, Euler's infinite product, Wendel's inequality-based limit, Raabe's formula, Binet's function, Stirling's formula, Gauss' multiplication formula, and Euler's constant. In Section 10, we provide a long list of properties of the polygamma function $\psi_{-2}$ that we can derive straightforwardly from this theory. Finally, in Section 11, we give some concluding remarks.

We observe that some alternative improvements of Bohr-Mollerup's theorem, in which both the convexity property and the asymptotic condition are somewhat relaxed, have been published in recent years. For a rather comprehensive list of references, see \cite[Section 3.3]{MarZen22}. Since this paper is a tutorial focused on the book \cite{MarZen22} rather than a survey paper, we will not elaborate here on these special improvements.

\section{Principal indefinite sums}

In this section we introduce the map, denoted $\Sigma$, that carries any function $g\colon\R_+\to\R$ satisfying the assumptions of Theorem~\ref{thm:int2} for some integer $p\geq 0$ into the unique function $f\colon\R_+\to\R$ defined in identity \eqref{eq:int2} and such that $f(1)=0$. For more details, see \cite[Chapter 5]{MarZen22}.

For any integer $p\geq 0$, we let $\cD^p$ denote the set of functions $g\colon\R_+\to\R$ that have the asymptotic property that the sequence $n\mapsto\Delta^pg(n)$ converges to zero. We also let $\cK^p$ denote the set of functions $g\colon\R_+\to\R$ that are eventually $p$-convex or eventually $p$-concave. We can easily show \cite[Chapter 4]{MarZen22} that $\cD^p\subset\cD^{p+1}$ and that $\cK^p\supset\cK^{p+1}$.

\begin{definition}\label{de:Sig4}
Let the map $\Sigma\colon\mathrm{dom}(\Sigma)\to\mathrm{ran}(\Sigma)$, with
$$
\mathrm{dom}(\Sigma) ~=~ \bigcup_{p\geq 0}(\cD^p\cap\cK^p),
$$
be defined by the condition
$$
g\in\cD^p\cap\cK^p \quad\Rightarrow\quad \Sigma g(x) ~=~ \lim_{n\to\infty}f^p_n[g](x),
$$
where $\mathrm{dom}(\Sigma)$ and $\mathrm{ran}(\Sigma)$ denote the domain and range of $\Sigma$, respectively.
\end{definition}


We observe that the map $\Sigma$ is well defined; indeed, if $g$ lies in both sets $\cD^p\cap\cK^p$ and $\cD^{p+1}\cap\cK^{p+1}$ for some integer $p\geq 0$, then necessarily
$$
\lim_{n\to\infty}f^{p+1}_n[g](x) ~=~ \lim_{n\to\infty}f^p_n[g](x),\qquad x>0.
$$

We also readily observe that the map $\Sigma$ is actually a bijection and its inverse is the restriction to $\mathrm{ran}(\Sigma)$ of the difference operator $\Delta$.

We can also show that
$$
\mathrm{ran}(\Sigma) ~=~ \bigcup_{p\geq 0}\{f\in\cK^p:\Delta f\in\cD^p\cap\cK^p~\text{and}~f(1)=0\}.
$$

Interestingly, Theorem~\ref{thm:int1} immediately provides the following characterization result.

\medskip

\begin{center}
\begin{minipage}{0.8\linewidth}
\emph{If $f\colon\R_+\to\R$ is a solution to the equation $\Delta f=g$, then it is eventually $p$-convex or eventually $p$-concave if and only if $f=c+\Sigma g$ for some $c\in\R$.}
\end{minipage}
\end{center}

\smallskip

\begin{definition}
We say that the \emph{principal indefinite sum} of a function $g$ lying in $\mathrm{dom}(\Sigma)$ is the class of functions $c+\Sigma g$, where $c\in\R$.
\end{definition}

\begin{example}[The log-gamma function]\label{ex:PIS43LOGG}
If $g(x)=\ln x$, then we have $\Sigma g(x)=\ln\Gamma(x)$, and we simply write
$$
\Sigma\ln x ~=~ \ln\Gamma(x),\qquad x>0.
$$
Thus, the principal indefinite sum of the function $x\mapsto \ln x$ is the class of functions $x\mapsto c+\ln\Gamma(x)$, where $c\in\R$. With some abuse of language, we can say that the principal indefinite sum of the log function is the log-gamma function.
\end{example}

\begin{example}[The polygamma function $\psi_{-2}$]\label{ex:PIS43LOGG22}
The function $g\colon\R_+\to\R$ defined by the equation
$$
g(x) ~=~ x\ln x -x+\frac{1}{2}\ln(2\pi),\qquad x>0,
$$
clearly lies in $\cD^2\cap\cK^2$. Its principal indefinite sum is the eventually $2$-convex function
$$
\Sigma g(x) ~=~ \psi_{-2}(x)-\psi_{-2}(1),\qquad x>0,
$$
where $\psi_{-2}$ is the polygamma function defined in Example~\ref{ex:PolyG}.
\end{example}

\section{Analogue of Gauss' limit and Eulerian form}

If the function $g\colon\R_+\to\R$ lies in $\cD^p\cap\cK^p$ for some integer $p\geq 0$, then by Definition~\ref{de:Sig4} we have
\begin{equation}\label{eq:GaussLim62}
\Sigma g(x) ~=~ \lim_{n\to\infty} f^p_n[g](x),\qquad x>0.
\end{equation}
As already discussed above, this latter identity is precisely the analogue of Gauss' limit for the gamma function. Moreover, it can be proved that the sequence $n\mapsto f^p_n[g]$ converges uniformly on any bounded subset of $\R_+$ to $\Sigma g$ (see Theorem~\ref{thm:int2}).

More generally, it was shown \cite[Section 7.1]{MarZen22} that, if $g$ is $r$ times continuously differentiable and lies in $\cD^p\cap\cK^{\max\{p,r\}}$ for some integer $r\geq 0$, then $\Sigma g$ is $r$ times continuously differentiable and the sequence $n\mapsto D^rf^p_n[g]$ converges uniformly on any bounded subset of $\R_+$ to $D^r\Sigma g$. In particular, both sides of \eqref{eq:GaussLim62} can be differentiated up to $r$ times and we have
$$
D^r\Sigma g(x) ~=~ \lim_{n\to\infty}D^rf_n^p[g](x),\qquad x>0.
$$

Moreover, if $g$ is continuous, then the function $f_n^p[g](x)-\Sigma g(x)$ can be integrated on any bounded interval of $[0,\infty)$ and the integral converges to zero as $n\to\infty$ (see \cite[Section 5.3]{MarZen22}).

Interestingly, the limit in \eqref{eq:GaussLim62} can be equivalently written in a series form. For instance, when $g(x)=\ln x$, the series representation of $\Sigma g$, once converted into the multiplicative notation, is precisely the following \emph{Euler product form} of the gamma function
\begin{equation}\label{eq:EulF}
\Gamma(x) ~=~ \frac{1}{x}\,\prod_{n=1}^{\infty}\frac{(1+1/n)^x}{1+x/n}{\,},\qquad x>0.
\end{equation}
This general observation is stated in the next theorem \cite[Section 8.1]{MarZen22}, which also shows that, under suitable assumptions, the series can be integrated and differentiated term by term.

\begin{theorem}[Eulerian form]\label{thm:Eul554}
Let $g$ lie in $\cD^p\cap\cK^p$ for some integer $p\geq 0$. Then the following assertions hold.
\begin{enumerate}
\item[(a)] For any $x>0$ we have
$$
\Sigma g(x) ~=~ -g(x)+\sum_{j=1}^p\tchoose{x}{j}{\,}\Delta^{j-1}g(1) - \sum_{n=1}^{\infty}\left(g(x+n)-\sum_{j=0}^p\tchoose{x}{j}{\,}\Delta^jg(n)\right)
$$
and the series converges uniformly on any bounded subset of $[0,\infty)$.
\item[(b)] If $g$ is continuous, then so is $\Sigma g$ and the series above can be integrated term by term on any bounded interval of $[0,\infty)$.
\item[(c)] If $g$ is $r$ times continuously differentiable and lies in $\cK^{\max\{p,r\}}$ for some integer $r\geq 0$, then $\Sigma g$ is $r$ times continuously differentiable and the series above can be differentiated term by term up to $r$ times.
\end{enumerate}
\end{theorem}

\begin{example}[The log-gamma function]
Consider the functions $g(x)=\ln x$ and $\Sigma g(x)=\ln\Gamma(x)$ given in Example~\ref{ex:PIS43LOGG} with the value $p=1$. Then, identity \eqref{eq:GaussLim62} clearly reduces to the additive version of Gauss' limit \eqref{eq:GaussLim6689}, that is
\begin{equation}\label{eq:GauLAddV}
\ln\Gamma(x) ~=~ \lim_{n\to\infty}\left(\ln(n-1)!+x\ln n-\sum_{k=0}^{n-1}\ln(x+k)\right).
\end{equation}
Similarly, using Theorem~\ref{thm:Eul554} we retrieve the additive version of the Euler product form \eqref{eq:EulF} of the gamma function, that is
$$
\ln\Gamma(x) ~=~ -\ln x-\sum_{n=1}^{\infty}\left(\ln(x+n)-\ln n-x\ln\left(1+\frac{1}{n}\right)\right).
$$
Moreover, the convergence is uniform on any bounded subset of $\R_+$.
\end{example}

\begin{example}[The polygamma function $\psi_{-2}$]\label{ex:PolySti4}
Consider the functions $g$ and $\Sigma g$ given in Example~\ref{ex:PIS43LOGG22} with the value $p=2$. We first observe that
$$
\Delta g(n) - \ln n ~\to ~ 0\qquad\text{as $n\to\infty$},
$$
and hence \eqref{eq:GaussLim62} yields the following identity
\begin{eqnarray*}
\psi_{-2}(x) &=& \lim_{n\to\infty}\Bigg(-x\,\sum_{k=1}^{n-1}\ln k-\sum_{k=1}^{n-1}(x+k)\ln\left(1+\frac{x}{k}\right)- x\ln x \\
&& \null + x\left(g(n)+n\right)+{x\choose 2}\ln n\Bigg),
\end{eqnarray*}
where the first sum clearly reduces to $\ln\Gamma(n)$. Now, using the additive version of Stirling's formula \eqref{eq:Stir46}, i.e.,
$$
\ln\Gamma(n)-g(n)+\frac{1}{2}\ln n ~\to ~ 0\qquad\text{as $n\to\infty$},
$$
we can easily see that the previous identity reduces to
$$
\psi_{-2}(x) ~=~ \lim_{n\to\infty}\left(nx-x\ln x+(\ln n)\frac{x^2}{2}-\sum_{k=1}^{n-1}(x+k)\ln\left(1+\frac{x}{k}\right)\right).
$$
Thus, this latter identity is a simplified form of the analogue of Gauss' limit for the gamma function. Interestingly, it can also be obtained directly by integrating both sides of \eqref{eq:GauLAddV}. The corresponding Eulerian form can be obtained similarly; we get
$$
\psi_{-2}(x) ~=~ x-x\ln x+\sum_{n=1}^{\infty}\left(x+\frac{x^2}{2}\ln\left(1+\frac{1}{n}\right)-(x+n)\ln\left(1+\frac{x}{n}\right)\right).
$$
Moreover, the convergence is uniform on any bounded subset of $\R_+$.
\end{example}

\section{The generalized Wendel's inequality-based limit}

For any integer $p\geq 0$, any real number $a>0$, and any function $g\colon\R_+\to\R$, we define the function $\rho^p_a[g]\colon [0,\infty)\to\R$ by the equation\label{p:rho}
\begin{equation}\label{eq:deflambdapt}
\rho^p_a[g](x) ~=~ g(x+a)-\sum_{j=0}^{p-1}\tchoose{x}{j}\,\Delta^jg(a),\qquad x\geq 0,
\end{equation}
or equivalently,
$$
\rho^p_a[g](x) ~=~ g(x+a)-P_{p-1}[g](a,a+1,\ldots,a+p-1;x+a),\qquad x\geq 0,
$$
where the function
$$
x~\mapsto ~ P_{p-1}[g](a,a+1,\ldots,a+p-1;x)
$$
denotes the unique interpolating polynomial of $g$ with nodes at the $p$ points $a,a+1,\ldots,a+p-1$. Thus, the quantity $\rho^p_a[g](x)$ is precisely the corresponding interpolation error at $x+a$.

We now present an important convergence result, which was established in \cite[Section 6.1]{MarZen22}.

\begin{theorem}[Generalized Wendel's inequality-based limit]\label{thm:Wendel9487}
For any integer $p\geq 0$, any real number $a\geq 0$, and any function $g$ lying in $\cD^p\cap\cK^p$, we have
$$
\rho^{p+1}_x[\Sigma g](a) ~\to ~0\qquad\text{as $x\to\infty$}{\,},
$$
or equivalently,
$$
\Sigma g(x+a)-\Sigma g(x)-\sum_{j=1}^p\tchoose{a}{j}\,\Delta^{j-1}g(x) ~\to ~0\qquad\text{as $x\to\infty$}{\,}.
$$
Moreover, if $g$ is $r$ times continuously differentiable and lies in $\cK^{\max\{p,r\}}$ for some integer $r\geq 0$, then this convergence result still holds if we differentiate the left-hand side with respect to $x$ up to $r$ times.
\end{theorem}

Applying Theorem~\ref{thm:Wendel9487} to the functions $g(x)=\ln x$ and $\Sigma g(x)=\ln\Gamma(x)$, with $p=1$, we immediately obtain
\begin{equation}\label{eq:6Rho2Ln0G24}
\rho_x^2[\Sigma\ln](a) ~=~ \ln\Gamma(x+a)-\ln\Gamma(x)-a\ln x
\end{equation}
and hence also the following well-known limit for any $a\geq 0$ (see, e.g., Titchmarsh~\cite{Tit39})
$$
\lim_{x\to\infty}\frac{\Gamma(x+a)}{\Gamma(x){\,}x^a} ~=~ 1{\,}.
$$
This latter result was also proved by Wendel~\cite{Wen48}, who first provided the following double inequality
$$
\left(1+\frac{a}{x}\right)^{a-1} \leq ~ \frac{\Gamma(x+a)}{\Gamma(x){\,}x^a} ~\leq ~ 1{\,}, \qquad x>0{\,},\quad 0\leq a\leq 1{\,},
$$
or equivalently, in the additive notation,
$$
(a-1)\ln\left(1+\frac{a}{x}\right) ~\leq ~ \rho_x^2[\Sigma\ln](a) ~\leq ~ 0{\,}, \qquad x>0{\,},\quad 0\leq a\leq 1{\,},
$$
which explains the name given to Theorem~\ref{thm:Wendel9487}.

\begin{example}[The polygamma function $\psi_{-2}$]
Let us apply Theorem~\ref{thm:Wendel9487} to the functions $g$ and $\Sigma g$ given in Example~\ref{ex:PIS43LOGG22} with the value $p=2$. Observing that
$$
\Delta g(x) - \ln x ~\to ~ 0\qquad\text{as $x\to\infty$},
$$
we can easily obtain the following limit for any $a\geq 0$,
$$
\psi_{-2}(x+a)-\psi_{-2}(x)-a\left(g(x)-\frac{1}{2}\ln x\right)-\frac{a^2}{2}\ln x ~\to ~ 0\qquad\text{as $x\to\infty$}.
$$
Using Stirling's formula as in Example~\ref{ex:PolySti4}, we finally obtain
$$
\psi_{-2}(x+a)-\psi_{-2}(x)-a\ln\Gamma(x)-\frac{a^2}{2}\ln x ~\to ~ 0\qquad\text{as $x\to\infty$},
$$
or equivalently,
$$
\int_0^a\ln\Gamma(x+t){\,}dt -a\ln\Gamma(x)-\frac{a^2}{2}\ln x ~\to ~ 0\qquad\text{as $x\to\infty$}.\qedhere
$$
\end{example}

\section{Analogue of Raabe's formula}

We now introduce a generalization of Raabe's formula. To this end, we first define the \emph{asymptotic constant} $\sigma[g]$ associated with a continuous function $g\colon\R_+\to\R$ lying in $\mathrm{dom}(\Sigma)$ by the equation
$$
\sigma[g] ~=~ \int_1^2\Sigma g(t){\,}dt ~=~ \int_0^1\Sigma g(t+1){\,}dt.
$$
For background on this concept, see \cite[Section 6.2]{MarZen22}.

Using this definition, we can immediately derive the following identity
\begin{equation}\label{eq:AnaRaa}
\int_x^{x+1}\Sigma g(t){\,}dt ~=~ \sigma[g]+\int_1^xg(t){\,}dt,\qquad x>0.
\end{equation}
Indeed, both sides of this identity are functions of $x$ that have the same derivative and the same value at $x=1$.

For instance, when $g(x)=\ln x$, we obtain
\begin{eqnarray*}
\sigma[g] &=& \int_0^1\ln\Gamma(t+1){\,}dt ~=~ \int_0^1(\ln\Gamma(t)+\ln t){\,}dt\\
&=& -1+\frac{1}{2}\ln(2\pi).
\end{eqnarray*}
Moreover, combining this value with \eqref{eq:AnaRaa} we obtain the following more general identity
\begin{equation}\label{eq:Raa55}
\int_x^{x+1}\ln\Gamma(t){\,}dt ~=~ x\ln x-x+\frac{1}{2}\ln(2\pi).
\end{equation}
This latter identity is known by the name \emph{Raabe's formula} (see, e.g., \cite[p.~46]{Rem98}). Thus, identity \eqref{eq:AnaRaa} provides for the function $\Sigma g$ an analogue of Raabe's formula.

\begin{example}[The polygamma function $\psi_{-2}$]\label{ex:51rrt}
Let us consider the functions $g$ and $\Sigma g$ given in Example~\ref{ex:PIS43LOGG22}. One can show \cite[p.~196]{Ada98} that
$$
\int_0^1\psi_{-2}(t){\,}dt ~=~ \ln A+\frac{1}{4}\ln(2\pi),
$$
where $A$ is the Glaisher-Kinkelin constant (see, e.g., Finch \cite[Section 2.15]{Fin03}). From this identity we immediately derive (see also \cite[Section 10.3]{MarZen22})
\begin{eqnarray*}
\sigma[g] &=& \int_0^1\psi_{-2}(t+1){\,}dt - \psi_{-2}(1) ~=~ \int_0^1(\psi_{-2}(t)+g(t)){\,}dt - \psi_{-2}(1)\\
&=& \ln A+\frac{1}{4}\ln(2\pi)-\frac{3}{4}{\,}.
\end{eqnarray*}
Identity \eqref{eq:AnaRaa} then provides for the function $\Sigma g$ the following analogue of Raabe's formula
$$
\int_x^{x+1}\psi_{-2}(t){\,}dt ~=~ \psi_{-2}(1)+\sigma[g]+\int_1^x g(t){\,}dt,
$$
or equivalently,
$$
\int_x^{x+1}\psi_{-2}(t){\,}dt ~=~ \frac{1}{2}{\,}x^2\ln x-\frac{3}{4}{\,}x^2+\frac{1}{4}\left(2x+1\right)\ln(2\pi)+\ln A.\qedhere
$$
\end{example}

\section{Generalized Binet's function}

Recall that the \emph{Binet function} related to the log-gamma function is the function $J\colon\R_+\to\R$ defined by the equation (see, e.g., Henrici \cite[p.~39]{Hen77})
\begin{equation}\label{eq:Binet5The}
J(x) ~=~ \ln\Gamma(x)-\frac{1}{2}\ln(2\pi)+x-\left(x-\frac{1}{2}\right)\ln x,\qquad x>0.
\end{equation}
Using identity \eqref{eq:6Rho2Ln0G24} and Raabe's formula \eqref{eq:Raa55}, we can easily provide the following integral form of Binet's function
$$
J(x) ~=~ -\int_0^1 \rho_x^2[\Sigma\ln](t){\,}dt,\qquad x>0.
$$

This latter formula motivates the following definition, in which we introduce a very useful generalization of Binet's function \cite[Section 6.3]{MarZen22}.

\begin{definition}[Generalized Binet's function]
For any integer $p\geq 0$ and any continuous function $g\colon\R_+\to\R$ lying in $\cD^p\cap\cK^p$, we define the function
$$
J^{p+1}[\Sigma g]\colon\R_+\to\R
$$
by the equation
\begin{equation}\label{eq:Binet643780}
J^{p+1}[\Sigma g](x) ~=~ -\int_0^1\rho_x^{p+1}[\Sigma g](t){\,}dt,\qquad x>0.
\end{equation}
We say that the function $J^{p+1}[\Sigma g](x)$ is the \emph{generalized Binet function} associated with the function $\Sigma g$ and the parameter $p+1$.
\end{definition}

Taking $g(x)=\ln x$ and $p=1$ in identity \eqref{eq:Binet643780}, we simply retrieve the Binet function
$$
J(x) ~=~ J^2[\Sigma\ln](x)
$$
related to the log-gamma function, as defined in \eqref{eq:Binet5The}.

Now, combining \eqref{eq:deflambdapt} with \eqref{eq:AnaRaa} and \eqref{eq:Binet643780}, we easily obtain the following explicit form of the generalized Binet function:
\begin{equation}\label{eq:BinSg68}
J^{p+1}[\Sigma g](x) ~=~ \Sigma g(x)-\sigma[g]-\int_1^xg(t){\,}dt + \sum_{j=1}^pG_j\Delta^{j-1} g(x),\qquad x>0,
\end{equation}
where $G_j$ is the $j$th \emph{Gregory coefficient} \cite{Bla16,MerSprVer06} defined by
$$
G_j ~=~ \int_0^1\tchoose{t}{j}{\,}dt.
$$

\begin{example}[The polygamma function $\psi_{-2}$]\label{ex:62rrt}
Consider the functions $g$ and $\Sigma g$ given in Example~\ref{ex:PIS43LOGG22} with the value $p=2$. Using identity \eqref{eq:BinSg68}, we obtain the following generalized Binet function
\begin{eqnarray*}
J^3[\Sigma g](x) &=& \psi_{-2}(x)-\frac{1}{12}{\,}(x+1)\ln(x+1)+\frac{1}{12}{\,}(3x-1)^2\\
&& \null -\frac{1}{12}{\,}x(6x-7)\ln x-\frac{1}{2}{\,}x\ln(2\pi)-\ln A.\qedhere
\end{eqnarray*}
\end{example}

\section{Generalized Stirling's formula}

We observe that the Binet function\index{Binet's function} $J(x)=J^2[\Sigma\ln](x)$ defined in \eqref{eq:Binet5The} clearly satisfies the following identity
$$
\Gamma(x) ~=~ \sqrt{2\pi}{\,}x^{x-\frac{1}{2}}{\,}e^{-x+J(x)}
$$
and hence Stirling's formula \eqref{eq:Stir46} simply states that $J(x)\to 0$ as $x\to\infty$. This observation is at the root of the following generalization of Stirling's formula \cite[Section 6.4]{MarZen22}.

\begin{theorem}[Generalized Stirling's formula]\label{thm:dgf7dds}
For any integer $p\geq 0$ and any continuous function $g\colon\R_+\to\R$ lying in $\cD^p\cap\cK^p$, the function $J^{p+1}[\Sigma g]$ vanishes at infinity. That is, using \eqref{eq:BinSg68},
$$
\Sigma g(x) -\int_1^x g(t){\,}dt +\sum_{j=1}^pG_j\Delta^{j-1}g(x) ~\to ~ \sigma[g]\qquad\text{as $x\to\infty$}.
$$
Moreover, if $g$ is $r$ times continuously differentiable and lies in $\cK^{\max\{p,r\}}$ for some integer $r\geq 0$, then this convergence result still holds if we differentiate both sides with respect to $x$ up to $r$ times.
\end{theorem}

Thus stated, the generalized Stirling formula enables one to investigate the asymptotic behavior of the function $\Sigma g$ for large values of its argument. When $g(x)=\ln x$ and $p=1$, we immediately retrieve the original Stirling formula \eqref{eq:Stir46}.

\begin{example}[The polygamma function $\psi_{-2}$]
Consider the functions $g$ and $\Sigma g$ given in Example~\ref{ex:PIS43LOGG22} with the value $p=2$. The corresponding generalized Stirling formula states that the function $J^3[\Sigma g]$ given in Example~\ref{ex:62rrt} vanishes at infinity. Using the fact that
$$
(x+1)\ln(x+1)-(x+1)\ln x -1 ~\to ~ 0 \qquad\text{as $x\to\infty$},
$$
this result can be restated as follows
$$
\psi_{-2}(x) - \frac{1}{12}(6x^2-6x+1)\ln x +\frac{1}{4}(3x-2)x-\frac{1}{2}{\,}x\ln(2\pi) ~\to ~ \ln A
$$
as $x\to\infty$. Differentiating both sides of this convergence result, we immediately retrieve the original Stirling formula.
\end{example}

\section{Analogue of Gauss' multiplication formula}

The additive version of Gauss' multiplication formula \eqref{eq:GaussMult} can be stated as follows; for any integer $m\geq 1$ we have
$$
\sum_{j=0}^{m-1}\ln\Gamma\left(\frac{x+j}{m}\right) ~=~ \ln\Gamma(x)-\left(x-\frac{1}{2}\right)\ln m+\frac{m-1}{2}\ln(2\pi),\qquad x>0.
$$

A generalization of this formula exists for any continuous function $\Sigma g$ lying in $\ran(\Sigma)$. It is stated in the following theorem \cite[Section 8.6]{MarZen22}.

\begin{theorem}[Analogue of Gauss' multiplication formula]\label{thm:AnGauMul57}
Let $m\geq 1$ be an integer and let $g\colon\R_+\to\R$ be a continuous function lying in $\mathrm{dom}(\Sigma)$. Define also the function $g_m\colon\R_+\to\R$ by the equation
$$
g_m(x) ~=~ g\left(\frac{x}{m}\right),\qquad x>0.
$$
Then the function $g_m$ also lies in $\mathrm{dom}(\Sigma)$. Moreover, for any $x>0$ we have
$$
\sum_{j=0}^{m-1}\Sigma g\left(\frac{x+j}{m}\right) ~=~ \Sigma g_m(x)+m{\,}\sigma[g]-\sigma[g_m]-\,\int_1^mg_m(t){\,}dt.
$$
\end{theorem}

Applying Theorem~\ref{thm:AnGauMul57} to the function $g(x)=\ln x$, we retrieve the original Gauss multiplication formula in its additive version. Let us now consider the case of the polygamma function $\psi_{-2}$.

\begin{example}[The polygamma function $\psi_{-2}$]\label{ex:82sad}
Let us apply Theorem~\ref{thm:AnGauMul57} to the functions $g$ and $\Sigma g$ given in Example~\ref{ex:PIS43LOGG22} with the value $p=2$. For any integer $m\geq 1$, we have
$$
g_m(x) ~=~ \frac{1}{m}{\,}g(x)-x\,\frac{\ln m}{m}+\frac{m-1}{2m}\ln(2\pi)
$$
and hence
$$
\Sigma g_m(x) ~=~ \frac{1}{m}{\,}\psi_{-2}(x)-{x\choose 2}\,\frac{\ln m}{m}+\frac{1}{2}\left(\frac{m-1}{m}{\,}x-1\right)\ln(2\pi).
$$
Using Theorem~\ref{thm:AnGauMul57}, after some algebra we obtain the following analogue of Gauss' multiplication formula
\begin{eqnarray*}
\sum_{j=0}^{m-1}\psi_{-2}\left(\frac{x+j}{m}\right) &=& \frac{1}{m}\,\psi_{-2}(x)-\frac{1}{12m}{\,}(6x^2-6x+1)\ln m \\
&& \null + (m-1)\ln(2\pi)\left(\frac{x}{2m}+\frac{1}{4}\right)+\left(m-\frac{1}{m}\right)\ln A.
\end{eqnarray*}
In particular, setting $x=1$ in this identity we obtain
$$
\sum_{j=1}^m\psi_{-2}\left(\frac{j}{m}\right) ~=~ \frac{1}{4}{\,}(m+1)\ln(2\pi)-\frac{1}{12m}\ln m+\left(m-\frac{1}{m}\right)\ln A.\qedhere
$$
\end{example}

\section{Generalized Euler's constant}

Recall that \emph{Euler's constant} (also called \emph{Euler-Mascheroni constant}) is defined as the limit
$$
\gamma ~=~ \lim_{n\to\infty}\left(\sum_{k=1}^n\frac{1}{k}-\ln n\right).
$$
This value actually represents the remainder in the numerical integration of the function $g(x)=1/x$ on the interval $[1,\infty)$ using the left rectangle method with the integer nodes $k=1,2,3,\ldots$ (see, e.g., Apostol~\cite{Apo99}).

A generalization of this value to any continuous function $g\colon\R_+\to\R$ lying in $\mathrm{dom}(\Sigma)$ was introduced in \cite[Section 6.8]{MarZen22} as follows.

\begin{definition}[Generalized Euler's constant]
Let $p\geq 0$ be an integer and let $g\colon\R_+\to\R$ be a continuous function lying in $\cD^p\cap\cK^p$. If $p\geq 1$, we also assume that $g$ does not lie in $\cD^{p-1}$. The \emph{generalized Euler constant} associated with the function $g$ is the number
$$
\gamma[g] ~=~ -J^{p+1}[\Sigma g](1),
$$
or equivalently, using \eqref{eq:BinSg68},
$$
\gamma[g] ~=~ \sigma[g]-\sum_{j=1}^p G_j\,\Delta^{j-1}g(1).
$$
\end{definition}


This definition can be justified by the following geometric interpretation. We can prove \cite[Section 6.8]{MarZen22} that
\begin{equation}\label{GeEuCoI5}
\gamma[g] ~=~ \int_1^{\infty}\left(\overline{P}_p[g](t)-g(t)\right)dt,
\end{equation}
where $\overline{P}_p[g]\colon [1,\infty)\to\R$ denotes the piecewise polynomial function whose restriction to the interval $[k,k+1)$, for any integer $k\geq 1$, is the interpolating polynomial of $g$ with nodes at $k,k+1,\ldots,k+p$; that is,
$$
\overline{P}_p[g](x) ~=~ P_p[g](k,k+1,\ldots,k+p;x),\qquad x\in [k,k+1).
$$
Moreover, if $g$ is $p$-convex or $p$-concave on $[1,\infty)$, then the graph of $g$ always lies over (or always lies under) that of $\overline{P}_p[g]$; and identity \eqref{GeEuCoI5} simply tells us that $|\gamma[g]|$ is the surface area between the two graphs on $[1,\infty)$.

\begin{example}
If $g(x)=\ln x$ and $p=1$, then we obtain
$$
\gamma[g] ~=~ \sigma[g] ~=~ -1+\frac{1}{2}\ln(2\pi) ~\approx ~ -0.081.
$$
The function $g$ is $1$-concave and its graph on $[1,\infty)$ always lies over that of the polygonal line $\overline{P}_1[g]$. The surface area between the two graphs is precisely $|\gamma[g]|$.
\end{example}

\begin{example}\label{ex:93rrt}
If
$$
g(x) ~=~ x\ln x -x+\frac{1}{2}\ln(2\pi)
$$
and $p=2$, then we obtain
$$
\gamma[g] ~=~ \ln A+\frac{1}{6}\ln 2-\frac{1}{3} ~\approx ~ 0.031.
$$
The function $g$ is $2$-concave and its graph on $[1,\infty)$ always lies under that of $\overline{P}_2[g]$. The surface area between the two graphs is precisely $\gamma[g]$.
\end{example}

\section{The polygamma function $\psi_{-2}$}

In the previous sections we have stated only some of our main results, starting from the generalization of Bohr-Mollerup's theorem and the analogue of Gauss' limit, and we have illustrated these results using our guiding example, the function $\psi_{-2}$. As mentioned in the introduction, many other results were established and illustrated in the book \cite{MarZen22}, where it was also demonstrated through several examples that all those results actually constitute a very powerful toolbox for the investigation of a large variety of functions.

To give the reader a taste of the scope of this new theory, in this section we simply present (without the detailed computations) what we can learn from it about the polygamma function $\psi_{-2}$.

Recall first that the polygamma function $\psi_{-2}\colon\R_+\to\R$ is defined by the equation (see Example~\ref{ex:PolyG})
$$
\psi_{-2}(x) ~=~ \int_0^x\ln\Gamma(t){\,}dt,\qquad x>0.
$$
Moreover, we have the identity (see Example~\ref{ex:PIS43LOGG22})
$$
\Sigma g(x) ~=~ \psi_{-2}(x)-\psi_{-2}(1),\qquad x>0,
$$
where $\psi_{-2}(1)=\frac{1}{2}\ln(2\pi)$ and
$$
g(x) ~=~ \Delta\psi_{-2}(x) ~=~ x\ln x-x+\psi_{-2}(1).
$$

Clearly, the function $g$ lies in $\cD^2$ and the function $\Sigma g$ lies in $\cD^3$. Moreover, both functions are infinitely many times differentiable. Furthermore, we can show \cite[Proposition 5.11]{MarZen22} that, for any integer $q\geq 1$, the function $g$ is eventually $(2q)$-concave and eventually $(2q+1)$-convex, while the function $\Sigma g$ is eventually $(2q)$-convex and eventually $(2q+1)$-concave.

\begin{remark}
We also observe that the function $\psi_{-2}$ is strongly related to the so-called \emph{hyperfactorial function} (or \emph{$K$-function}). Indeed, the latter is the function $K\colon\R_+\to\R_+$ defined by the equations
\begin{eqnarray*}
\ln K(x) &=& \zeta'(-1,x)-\zeta'(-1)\\
&=& \tchoose{x}{2} +\psi_{-2}(x)-x\,\psi_{-2}(1)\qquad\text{for $x>0$},
\end{eqnarray*}
where $\zeta(s,x)$ denotes the Hurwitz zeta function and $\zeta'(s,x)$ denotes its derivative with respect to the variable $s$ (see, e.g., \cite[p.~196]{Ada98} and \cite[Section 12.5]{MarZen22}). Thus defined, the hyperfactorial function clearly satisfies the identity $\Delta\ln K(x)=x\ln x$ on $\R_+$ (or equivalently, $K(x+1)=x^x{\,}K(x)$ on $\R_+$). Moreover, the analogue of Bohr-Mollerup's theorem states that the function $f(x)=\ln K(x)$ is the unique (up to an additive constant) eventually $2$-convex solution to the equation $\Delta f(x)=x\ln x$ on $\R_+$.
\end{remark}

\subsection{Analogue of Bohr-Mollerup's theorem}

The function $\psi_{-2}$ can be characterized as follows (see Example~\ref{ex:PolyG}).

\begin{theorem}
A function $f\colon\R_+\to\R$ is a solution to the equation $\Delta f=g$ that lies in $\cK^2$ if and only if it is of the form $f=c+\psi_{-2}$, where $c\in\R$.
\end{theorem}
We also have the following alternative characterization \cite[Section 3.1]{MarZen22}.
\begin{theorem}
A function $f\colon\R_+\to\R$ is a solution to the equation $\Delta f=g$ that has the property that, for each $x>0$, the sequence
$$
n ~\mapsto ~ f(x+n)-f(n)-x{\,}\ln\Gamma(n)-\frac{x^2}{2}\ln n
$$
converges to zero, if and only if it is of the form $f=c+\psi_{-2}$, where $c\in\R$.
\end{theorem}

\subsection{Asymptotic constant and generalized Euler's constant}

We have the following values (see Examples~\ref{ex:51rrt} and \ref{ex:93rrt})
\begin{eqnarray*}
\sigma[g] &=& \ln A+\frac{1}{4}\ln(2\pi)-\frac{3}{4}{\,},\\
\gamma[g] &=& \ln A+\frac{1}{6}\ln 2-\frac{1}{3}{\,}.
\end{eqnarray*}
We also have the following integral representations \cite[Section 10.3]{MarZen22}
\begin{eqnarray*}
\sigma[g] &=& \frac{1}{2}{\,}g(1)-\frac{1}{2}\int_1^{\infty}\frac{B_2(\{t\})}{t}{\,}dt,\\
\gamma[g] &=& \int_1^{\infty}\left(-g(t)+g(\lfloor t\rfloor)+\frac{1}{2}\Delta g(\lfloor t\rfloor)-\frac{1}{12}\Delta^2g(\lfloor t\rfloor)\right)dt,
\end{eqnarray*}
where $\{t\}=t-\lfloor t\rfloor$ and $B_2$ is the Bernoulli polynomial $B_2(x)=x^2-x+1/6$.

\subsection{Analogue of Raabe's formula}

We have the following analogue of Raabe's formula (see Example~\ref{ex:51rrt})
$$
\int_x^{x+1}\psi_{-2}(t){\,}dt ~=~ \frac{1}{2}{\,}x^2\ln x-\frac{3}{4}{\,}x^2+\frac{1}{4}\left(2x+1\right)\ln(2\pi)+\ln A,\qquad x>0.
$$
Moreover, the function $f=\psi_{-2}$ is the unique continuous solution lying in $\cK^2$ to the equation (see \cite[Section 8.5]{MarZen22})
$$
\int_x^{x+1}f(t){\,}dt ~=~ \frac{1}{2}{\,}x^2\ln x-\frac{3}{4}{\,}x^2+\frac{1}{4}\left(2x+1\right)\ln(2\pi)+\ln A,\qquad x>0.
$$

\subsection{Generalized Binet's function}

We have the following generalized Binet function (see Example~\ref{ex:62rrt})
\begin{eqnarray*}
J^3[\Sigma g](x) &=& \psi_{-2}(x)-\frac{1}{12}{\,}(x+1)\ln(x+1)+\frac{1}{12}{\,}(3x-1)^2\\
&& \null -\frac{1}{12}{\,}x(6x-7)\ln x-\frac{1}{2}{\,}x\ln(2\pi)-\ln A.
\end{eqnarray*}

\subsection{Inequalities}

The following inequalities hold for any $x>0$ and any $a\geq 0$.
\begin{itemize}
\item\emph{Generalized Wendel's inequality} \cite[Section 6.1]{MarZen22}
\begin{eqnarray*}
0 &\leq & \mathrm{sign}(a(a-1)(a-2))\left(\psi_{-2}(x+a)-\psi_{-2}(x)-a{\,}g(x)-\tchoose{a}{2}\Delta g(x)\right)\\
&\leq & \left|\tchoose{a-1}{2}\right|(\Delta g(x+a)-\Delta g(x)) ~\leq ~ \lceil a\rceil\left|\tchoose{a-1}{2}\right|\Delta^2g(x).
\end{eqnarray*}
\item\emph{Generalized Webster's inequality} \cite[Appendix E]{MarZen22}
\begin{eqnarray*}
0 &\leq & \psi_{-2}(x+a+1)-\psi_{-2}(x+\lfloor a\rfloor +1)\\
&& \null -\{a\}{\,}g(x+\lfloor a\rfloor +1)-\tchoose{\{a\}}{2}\Delta g(x+\lfloor a\rfloor +1)\\
&\leq & \frac{1}{2}\{a\}\left(g(x+a)-g(x+\lfloor a\rfloor +1)-(\{a\}-1)\Delta g(x+\lfloor a\rfloor +1)\right),
\end{eqnarray*}
where $\{a\}=a-\lfloor a\rfloor$.
\item\emph{Generalized Gautschi's inequality} \cite[Section 10.3]{MarZen22}
\begin{eqnarray*}
(a-\lceil a\rceil)\,\ln\Gamma(x+\lceil a\rceil) &\leq & \psi_{-2}(x+a)-\psi_{-2}(x+\lceil a\rceil)\\
&\leq & (a-\lceil a\rceil){\,}g(x+\lfloor a\rfloor),
\end{eqnarray*}
provided $x+\lfloor a\rfloor\geq x_0$, where $x_0=1.461\ldots$ is the unique positive zero of the digamma function.
\item\emph{Generalized Stirling's formula-based inequality} \cite[Section 6.4]{MarZen22}
\begin{eqnarray*}
0 &\leq & -J^3[\Sigma g](x) ~\leq ~ \int_0^1\tchoose{t-1}{2}(\Delta g(x+t)-\Delta g(x)){\,}dt\\
&\leq & \frac{5}{12}\,\Delta^2g(x).
\end{eqnarray*}
\end{itemize}

We also have the following double inequality \cite[Appendix E]{MarZen22}
$$
\alpha(x) ~\leq ~ \psi_{-2}(x) ~\leq ~ \beta(x),\qquad x>0,
$$
where
\begin{eqnarray*}
\alpha(x) &=& \ln A-\frac{5}{18}+\frac{1}{24}{\,}x-\frac{5}{6}{\,}x^2+\frac{1}{2}{\,}x\ln(2\pi)-\frac{1}{12}{\,}x(x^2+12)\ln x\\
&& \null + \frac{1}{12}{\,}(x+1)(x^2+5x+1)\ln(x+1)
\end{eqnarray*}
and
\begin{eqnarray*}
\beta(x) &=& \ln A-\frac{1}{3}-\frac{3}{4}{\,}x^2+\frac{1}{2}{\,}x\ln(2\pi)-x\ln x\\
&& \null +\frac{1}{12}{\,}(x+1)(6x-1)\ln(x+1)+\frac{1}{12}{\,}(x+2)\ln(x+2).
\end{eqnarray*}
This double inequality provides a rather fine bracketing of the function $\psi_{-2}$ for large values of $x$. Indeed, we have
$$
\sup_{x\in\R_+}{\,}|\beta(x)-\alpha(x)| ~=~ \frac{1}{18}{\,}(3\ln 2-1) ~\approx ~ 0.06
$$
and
$$
\beta(x)-\alpha(x) ~=~ \frac{1}{16{\,}x}-\frac{13}{180{\,}x^2}+\frac{13}{144{\,}x^3}+O(x^{-4})\qquad\text{as $x\to\infty$}.
$$

\subsection{Generalized Stirling's and related formulas}

For any $a\geq 0$, we have the following limits as $x\to\infty$ \cite[Section 10.3]{MarZen22}
$$
\psi_{-2}(x+a)-\psi_{-2}(x)-a\ln\Gamma(x)-\frac{a^2}{2}\ln x ~\to ~ 0,
$$
$$
\psi_{-2}(x) - \frac{1}{12}(6x^2-6x+1)\ln x +\frac{1}{4}(3x-2)x-\frac{1}{2}{\,}x\ln(2\pi) ~\to ~ \ln A,
$$
$$
\psi_{-2}(x) -x\ln\Gamma(x)+\frac{1}{12}(6x^2-1)\ln x-\frac{1}{4}{\,}x(x+2) ~\to ~ \ln A-\frac{1}{12}{\,},
$$
$$
\frac{\psi_{-2}(x+a)}{x^2\ln x} ~\to ~ \frac{1}{2}{\,}.
$$

\subsection{Asymptotic expansions}

For any integers $m\geq 1$ and $q\geq 1$, we have the following expansion as $x\to\infty$ \cite[Section 10.3]{MarZen22}
\begin{eqnarray*}
\frac{1}{m}\sum_{j=0}^{m-1}\psi_{-2}\left(x+\frac{j}{m}\right) &=& \frac{1}{2}{\,}x^2\ln x-\frac{3}{4}{\,}x^2+\left(\frac{1}{2}{\,}x+\frac{1}{4}\right)\ln(2\pi)+\ln A\\
&& \null + \sum_{k=1}^q\frac{B_k}{m^kk!}{\,}g^{(k-1)}(x)+O(g^{(q)}(x)).
\end{eqnarray*}
Setting $m=1$ in this formula, we obtain
\begin{eqnarray*}
\psi_{-2}(x) &=& \frac{1}{2}{\,}x^2\ln x-\frac{3}{4}{\,}x^2+\left(\frac{1}{2}{\,}x+\frac{1}{4}\right)\ln(2\pi)+\ln A\\
&& \null + \sum_{k=1}^q\frac{B_k}{k!}{\,}g^{(k-1)}(x)+O(g^{(q)}(x)).
\end{eqnarray*}
For instance, we have
\begin{eqnarray*}
\psi_{-2}(x) &=& \frac{1}{12}{\,}(6x^2-6x+1)\ln x-\frac{1}{4}(3x-2)x+\frac{1}{2}{\,}x\ln(2\pi)+\ln A\\
&& \null + \frac{1}{720{\,}x^2}-\frac{1}{5040{\,}x^4}+\frac{1}{10080{\,}x^6} + O(x^{-8}).
\end{eqnarray*}

\subsection{Generalized Liu's formula}

For any $x>0$, we have \cite[Section 10.3]{MarZen22}
\begin{eqnarray*}
\psi_{-2}(x) &=& \frac{1}{12}{\,}(6x^2-6x+1)\ln x-\frac{1}{4}(3x-2)x+\frac{1}{2}{\,}x\ln(2\pi)+\ln A\\
&& \null + \frac{1}{2}\int_0^{\infty}\frac{B_2(\{t\})}{x+t}{\,}dt,
\end{eqnarray*}
where $\{t\}=t-\lfloor t\rfloor$ and $B_2$ is the Bernoulli polynomial $B_2(x)=x^2-x+1/6$.

\subsection{Limit, series, and integral representations}

We have the following formulas for any $x>0$ (see Example~\ref{ex:PolySti4} and \cite[Section 8.2]{MarZen22})
\begin{itemize}
\item\emph{Analogue of Gauss' limit}
$$
\psi_{-2}(x) ~=~ \lim_{n\to\infty}\left(nx-x\ln x+(\ln n)\frac{x^2}{2}-\sum_{k=1}^{n-1}(x+k)\ln\left(1+\frac{x}{k}\right)\right).
$$
\item\emph{Eulerian form}
$$
\psi_{-2}(x) ~=~ x-x\ln x+\sum_{n=1}^{\infty}\left(x+\frac{x^2}{2}\ln\left(1+\frac{1}{n}\right)-(x+n)\ln\left(1+\frac{x}{n}\right)\right).
$$
\item\emph{Weierstrassian form}
$$
\psi_{-2}(x) ~=~ -\gamma{\,}\frac{x^2}{2}+x-x\ln x+\sum_{n=1}^{\infty}\left(x+\frac{1}{2n}{\,}x^2-(x+n)\ln\left(1+\frac{x}{n}\right)\right),
$$
where $\gamma$ is Euler's constant.
\item\emph{Integral form}
$$
\psi_{-2}(x) ~=~ \int_0^x\ln\Gamma(t){\,}dt ~=~ x\ln\Gamma(x)-\int_0^x t\,\psi(t){\,}dt,
$$
where $\psi(x)=D\ln\Gamma(x)$ is the digamma function.
\end{itemize}

\subsection{Analogue of Gauss' multiplication formula}

For any $x>0$ and any integer $m\geq 1$, we have (see Example~\ref{ex:82sad})
\begin{eqnarray*}
\sum_{j=0}^{m-1}\psi_{-2}\left(\frac{x+j}{m}\right) &=& \frac{1}{m}\,\psi_{-2}(x)-\frac{1}{12m}{\,}(6x^2-6x+1)\ln m \\
&& \null + (m-1)\ln(2\pi)\left(\frac{x}{2m}+\frac{1}{4}\right)+\left(m-\frac{1}{m}\right)\ln A.
\end{eqnarray*}
Letting $x\to 0$ in this identity, we obtain
$$
\sum_{j=1}^{m-1}\psi_{-2}\left(\frac{j}{m}\right) ~=~ -\frac{1}{12m}\ln m+\frac{1}{4}{\,}(m-1)\ln(2\pi)+\left(m-\frac{1}{m}\right)\ln A.
$$
For instance, when $m=2$ we immediately derive the formula
$$
\psi_{-2}\left(\frac{1}{2}\right) ~=~ \frac{5}{24}\ln 2+\frac{1}{4}\ln\pi+\frac{3}{2}\ln A.
$$
Interestingly, we can also derive the following limit \cite[Section 10.3]{MarZen22}
$$
\lim_{m\to\infty}\left(\frac{1}{m^2}\,\psi_{-2}(mx)-\frac{1}{2}{\,}x^2\ln m\right) ~=~ \frac{1}{2}{\,}x^2\ln x-\frac{3}{4}{\,}x^2,\qquad x>0.
$$

\subsection{Gregory's formula-based series representation}

For any $x>0$, we have \cite[Section 10.3]{MarZen22}
$$
\psi_{-2}(x) ~=~ \frac{1}{2}{\,}x^2\ln x-\frac{3}{4}{\,}x^2+\left(\frac{1}{2}{\,}x+\frac{1}{4}\right)\ln(2\pi)+\ln A-\sum_{n=1}^{\infty}G_n\,\Delta^{n-1}g(x),
$$
where $G_n$ is the $n$th Gregory coefficient. Setting $x=1$ in this identity, we obtain the following analogue of \emph{Fontana-Mascheroni's series}
$$
\sum_{n=1}^{\infty}G_n\,\Delta^{n-1}g(1) ~=~ \sigma[g] ~=~ \ln A+\frac{1}{4}\ln(2\pi)-\frac{3}{4}{\,}.
$$

\subsection{Analogue of Wallis's product formula}

We have the following analogues of Wallis's product formula \cite[Section 10.3]{MarZen22}
\begin{eqnarray*}
\lim_{n\to\infty}\left(h_1(n)+\sum_{k=1}^{2n}(-1)^{k-1}g(k)\right) &=& \frac{1}{12}\ln 2-3\ln A,\\
\lim_{n\to\infty}\left(h_2(n)+\sum_{k=1}^{2n}(-1)^{k-1}\,\psi_{-2}(k)\right) &=& \ln A -\frac{1}{12}\ln 2,
\end{eqnarray*}
where
\begin{eqnarray*}
h_1(n) &=& \left(n+\frac{1}{4}\right)\ln n-n(1-\ln 2),\\
h_2(n) &=& n^2\ln(2n)-\frac{3}{2}{\,}n^2+\frac{1}{2}{\,}n\ln(2\pi)-\frac{1}{12}\ln n.
\end{eqnarray*}

\subsection{Generalized Webster's functional equation}

For any integer $m\geq 1$, there is a unique solution $f\colon\R_+\to\R$ to the equation
$$
\sum_{j=0}^{m-1}f\left(x+\frac{j}{m}\right) ~=~ g(x)
$$
that lies in $\cK^2$, namely \cite[Section 10.3]{MarZen22}
$$
f(x) ~=~ \psi_{-2}\left(x+\frac{1}{m}\right)-\psi_{-2}(x).
$$

\subsection{Analogue of Euler's series representation of $\gamma$}

The Taylor series expansion of $\psi_{-2}(x+1)$ about $x=0$ is \cite[Section 10.3]{MarZen22}
$$
\psi_{-2}(x+1) ~=~ \frac{1}{2}\ln(2\pi)-\gamma{\,}\frac{x^2}{2}+\sum_{n=3}^{\infty}(-1)^{n-1}\,\frac{\zeta(n-1)}{n(n-1)}{\,}x^n,\qquad |x|<1,
$$
where $z\mapsto\zeta(z)$ denotes the Riemann zeta function. Integrating both sides of this equation on $(0,1)$, we obtain the following analogue of \emph{Euler's series representation of $\gamma$}
$$
\sum_{n=2}^{\infty}(-1)^n\,\frac{\zeta(n)}{n(n+1)(n+2)} ~=~ \frac{1}{6}\,\gamma-\frac{3}{4}+\frac{1}{4}\ln(2\pi)+\ln A.
$$

\subsection{Analogue of Euler's reflection formula}

For any $x\in (0,1)$, we have \cite[Section 10.3]{MarZen22}
$$
\psi_{-2}(x)-\psi_{-2}(1-x) ~=~ x\ln\pi-\frac{1}{2}\ln(2\pi)-\int_0^x\ln\sin(\pi t){\,}dt.
$$

\section{Conclusion}

The authors have recently published an open access book \cite{MarZen22} that presents a significant generalization of Bohr-Mollerup's theorem to higher order convex functions. This generalization shows that a very rich spectrum of functions satisfy analogues of several classical properties of the gamma function, including Bohr-Mollerup's theorem itself, Euler's reflection formula, Gauss' multiplication theorem, Stirling's formula, and Weierstrass' canonical factorization.

In this tutorial paper, we have summarized the main results of this new theory and have illustrated these results as well as many others by applying them to the polygamma function $\psi_{-2}$ (i.e., the integral of the log-gamma function).

Actually, the uniqueness and existence theorems given in the introduction show that it is usually not very difficult to check whether a given function can be investigated through our results or not. If so, then we may gain a lot of insight into this function just by applying those results almost mechanically.

In writing this paper, our hope is to convince the reader that our theory offers a unifying approach that enables us to systematically handle a wide variety of functions all at once using elementary means.

Beyond this systematization aspect, this theory introduces some new important and useful objects. For instance, the map $\Sigma$ itself is a new concept that seems to be as fundamental as the basic antiderivative operation. Other concepts such as the asymptotic constant and the generalized Binet function also appear to be new fundamental objects that merit further investigation. These objects are used, e.g., in the remarkable generalized Stirling formula, but also in many other useful formulas such as the analogue of Raabe's formula and the analogue of Gauss' multiplication formula.

Lastly, this theory also revealed how natural and useful the higher order convexity properties are. Although these properties seem to be still rather poorly investigated in mathematical analysis, they clearly play a crucial role in this setting and hence also merit further study.

\section{Declarations}

\noindent\textbf{Ethical Approval} Not applicable.

\smallskip\noindent\textbf{Conflict of interest} None.

\smallskip\noindent\textbf{Authors' contributions} Both authors contributed equally to the manuscript and
revised its final form.

\smallskip\noindent\textbf{Funding} This research received no external funding.

\smallskip\noindent\textbf{Availability of data and materials} Not applicable.


\end{document}